\documentclass{dimacs-l}
\usepackage{amsmath}
\usepackage{amssymb}
\usepackage{amsthm}
\usepackage{epsf}

\let\cal\mathcal

\newcounter{AbcT}

\newtheorem {Theorem}    {Theorem}[section]
\newtheorem {Definition} {Definition} [section]
\newtheorem {Problem}    {Problem}
\newtheorem {Conjecture} {Conjecture}

\newtheorem {Proposition}[Theorem]    {Proposition}

\newtheorem {Example}    [Theorem]    {Example}

\newcommand {\proofs}     {\proof}

\let\QED\qed

\newcommand{\ignore}[1]{}
\newcommand{\br}{{\mbox{{\rm br}}}}
\newcommand{\rank}{{\mbox{{\rm rank}}}}

\def\pho{\rho}
\def\eps{\epsilon}
\def\E{{\bf{E}}}
\def\P{{\bf{P}}}

\def\M{{\bf{M}}}
\def\R{\hbox{I\kern-.2em\hbox{R}}}
\def\A{{\cal{A}}}

\def\F{{\cal{F}}}

\def\|{\, | \, }
\def\v0{{\bf 0}}

\def\0{\hat{0}}
\def\1{\hat{1}}

\def\path{{\tt path}}

\def\lam{\lambda}

\def\phi{\varphi}

\def\be{\begin{equation}}
\def\ee{\end{equation}}

\begin{document}

\title{Survey: Information flow on trees}
\author{Elchanan Mossel}
\address{Statistics, University of California, Berkeley, CA 94720}
\email{mossel@stat.berkeley.edu}

\begin{abstract}

Consider a tree network $T$, where each edge acts as an independent copy of
a given channel $M$, and information is propagated from the root.
For which $T$ and $M$ does the configuration
obtained at level $n$ of $T$ typically contain significant
information on the root variable?

This model appeared independently in biology, information theory, and
statistical physics.
Its analysis uses techniques
from the theory of finite markov chains, statistics, statistical physics,
information theory, cryptography and noisy computation.
In this paper, we survey developments and challenges related to this problem.

\end{abstract}

\copyrightinfo{2004}
  {American Mathematical Society}

\maketitle

\section{Introduction}
Consider a process on a tree in which information is transmitted from
the root of the tree to all the nodes of the tree.
Each node inherits information from its parent with some probability of error.
The transmission process is assumed to have identical distribution on all the
edges, and different edges of the tree are assumed to act independently.

As this process represents propagation of a genetic property from ancestor
to its descendants, it was studied in genetics, see e.g. \cite{Ca,SC}.
In communication theory, this process represents a
communication network on the tree where information is
transmitted from the root of the tree.
Earlier, the process was studied in statistical physics, see e.g.
\cite{Sp,Hi,BRZ}.

The basic question we address in this survey is:
Does the configuration
obtained at level $n$ of $T$ typically contain significant
information on the root variable?
The theory of finite markov chains
implies that if the underlying markov
chain is ergodic (i.e. irreducible and aperiodic),
then the variable at a single node at level $n$ and the variable at the
root are asymptotically independent as $n \to \infty$.
However, for the tree process, information is duplicated,
so it is conceivable that the configuration at level
$n$ contains a significant amount of information on the root variable.

In Section \ref{sec:definitions} a precise formulation of the problem is given.
In Section \ref{sec:ising} we discuss the problem for symmetric binary channels
(which correspond to Ising models with no external field).
This is the family of channels for which the most is known.
In particular, is subsection \ref{subsec:alg} we compare
various reconstruction algorithms for symmetric binary channels.

Suppose that instead of the {\em configuration} at level $n$, we are given the {\em census} of
the configuration at level $n$. In section \ref{sec:census} we see how the spectral properties of $M$
determine if the census  is asymptotically independent of the root.
In Section \ref{sec:potts} we present bounds for the problem for Potts models, while
in Section \ref{sec:strange} we present some examples of channels which are related to
secret sharing and demonstrate that level $n$ may contain significant information on the
root even if the census of the level contains no such information.

In Section \ref{sec:terminology} we discuss related problems and terminology.
Many unsolved problems are presented throughout the paper.

\medskip\noindent
{\bf Disclaimer:} This paper is a survey: most of the results are presented
without a proof; for others only a sketch is given.

\section{The reconstruction problem} \label{sec:definitions}
{\bf Definition of the process.}
Denote the underlying tree by $T = (V,E)$.
The information flow on each edge is given by
a channel on a finite alphabet $\A=\{1,\ldots,k\}$.
Let $\M_{i,j}$ be the transition probability from $i$ to $j$;
$M$ be the random function (or channel) which satisfies for all $i$
and $j$ that
$\P[M(i) = j] = \M_{i,j}$, and
$\lam_2(M)$ be the eigenvalue of $\M$ which has the
second largest absolute value ($\lam_2(M)$ is in general a complex number).
We assume throughout the paper that $M$
defines an ergodic markov chain (irreducible and aperiodic).
At the root $\rho$ one of the symbols of $\A$ is chosen according to
some initial distribution. 
We denote this (random) symbol by $\sigma_{\rho}$.
This symbol is then propagated in the tree in the following
way. Given that the parent of $v$, denoted $v'$, has value $\sigma_{v'}$,
the probability that $\sigma_v$ is $j$ is given by $\M_{\sigma_{v'},j}$.
More formally, for each vertex $v$ having as a parent $v'$, we let
$\sigma_v = M_{v',v}(\sigma_{v'})$, where the $\{M_{v',v}\}$
are independent copies of $M$.
Equivalently, for a vertex $v$, let $v'$ be the parent of $v$, and let
$\Gamma(v)$ be the set of all vertices which are connected to $\pho$
through paths which do not contain $v$. Then the process
satisfies:
\[
\P[\sigma_v = j | (\sigma_w)_{w \in \Gamma(v)}] = \P[\sigma_v = j | \sigma_{v'}]
 = \M_{\sigma_{v'},j}.
\]
Let $d(,)$ denote the graph-metric distance on $T$, and
$L_n = \{v \in V: d(\rho,v) = n\}$ be the $n$'th level of the tree.
For $v \in V$ and $e=(v,w) \in E$ we denote $|v| = d(\rho,v)$ and
$|e| = \max\{|v|,|w|\}$.
We denote by $\sigma_n = (\sigma(v))_{v \in L_n}$ the symbols at the
$n$'th level of the tree.
We let $c_n = (c_n(1),\ldots,c_n(k))$ where
\[
c_n(i) = \#\{v \in L_n : \sigma(v) = i\}.
\]
In other words, $c_n$ is the {\sl census} of the $n$'th level.
Note that both $(\sigma_n)_{n=1}^{\infty}$ and
$(c_n)_{n=1}^{\infty}$ are markov chains.

\noindent
{\bf Reconstruction solvability.}
For distributions $P$ and $Q$ on the same space, the total
variation distance between $P$ and $Q$ is
\begin{equation} \label{eq:totalvar}
D_V(P,Q) = \frac{1}{2} \sum_{\sigma} |P(\sigma) - Q(\sigma)|.
\end{equation}
\begin{Definition} \label{def:reconstruction}
The reconstruction problem for $T$ and $M$ is {\em solvable} if
there exist $i,j \in \A$ for which
\begin{equation} \label{cond:arb_l1}
\lim_{n \to \infty} D_V(\P_n^i,\P_n^j) > 0,
\end{equation}
where
$\P_n^\ell$ denotes the conditional distribution of
$\sigma_n$ given that $\sigma_{\rho} = \ell$.
We define {\em census} solvability similarly, where the measures
$\P_n^{\ell}$ are replaced by measures $\widetilde{\P}_n^{\ell}$,
which are conditional distributions of $c_n$ given that $\sigma_{\rho} =
\ell$.
\end{Definition}
A stronger definition than Definition \ref{def:reconstruction} is obtained
by replacing ``there exists $i,j$'' by ``for all $i,j$''. We choose the 
weaker definition as we are interested to know if some information is propagated from the root to the boundary, see also Proposition \ref{prop:equiv} below.    

\noindent
{\bf Equivalent definitions.}
If the reconstruction problem is solvable, then $\sigma_n$
contains significant information on the root variable.
This may be expressed in several equivalent ways.
Assume that the variable at the root, $\sigma_{\rho}$, is chosen according
to some initial distribution $(\pi_i)_{i \in \A}$, and let $\P^{\pi}$
denote the corresponding probability measure.
The maximum-likelihood algorithm, which is the
optimal reconstruction algorithm of $\sigma_{\rho}$ given $\sigma_n$,
is successful with probability
\begin{eqnarray} \nonumber
\Delta_n(\pi) &=& \sum_{\sigma}
\P^{\pi}[\sigma_n = \sigma]
\max_{i \in \A} \P^{\pi}[\sigma_{\rho} = i | \sigma_n = \sigma] \\ \nonumber
&\geq&
\max_{i \in \A} \sum_{\sigma} \P^{\pi}[\sigma_n = \sigma] \P^{\pi}[\sigma_{\rho} = i | \sigma_n = \sigma] =
\max_{i \in \A} \pi_i.
\end{eqnarray}
Note that it is possible to reconstruct $\sigma_{\rho}$
with probability $\max_i \pi_i$ even when $\sigma_n$ is unknown
(using the algorithm which always reconstructs the $i$ which maximizes
$\pi_i$).
It is therefore natural to consider $\Delta_n(\pi) - \max_i \pi_i$ as a measurement for the dependency
between $\sigma_n$ and $\sigma_{\rho}$.

Let $H$ be the entropy function, and
let $I(X,Y) = H(X) + H(Y) - H(X,Y)$ be the mutual-information operator
(see e.g. \cite{CT} for definitions and basic properties).

For a sequence of random variables $X_n$ defined on the same probability 
space, let $F_n$ be the $\sigma$-algebra defined by $(X_m)_{m \geq n}$, i.e., 
$F_n$ is the minimal $\sigma$-algebra such that all the variables 
$(X_m)_{m \geq n}$ are measurable with respect to $F_n$. 
Let $F_{\infty} = \cap_{n=1}^{\infty} F_n$.
We say that the sequence $X_n$ has a trivial tail, if all the measurable sets with respect to $F_{\infty}$ have probability either $0$ or $1$.
Otherwise, we say the the sequence has a non-trivial tail.

In the theory of markov random fields the notion of tail triviality is 
closely related to the extremality of the measure, see e.g. \cite{Ge}.
   
The following equivalence follows from the fact that $\sigma_n$ is
a markov chain (see e.g. \cite{M:lam2}):
\begin{Proposition} \label{prop:equiv}
Let $T$ be an infinite tree and $M$ a channel.
Then the following conditions are equivalent
(where $\pi$ denotes initial distribution for $\sigma_{\rho}$):
\begin{enumerate}
\item
The reconstruction problem is solvable
\item
There exists a $\pi$ for which
$\lim_{n \to \infty} I(\sigma_{\rho},\sigma_n) > 0$.
\item
If $\pi$ is the uniform distribution on $\A$, then
$\lim_{n \to \infty} I(\sigma_{\rho},{\sigma}_n) > 0$.
\item
For any distribution $\pi$ with $\min_i \pi_i > 0$, it holds that
$\lim_{n \to \infty} I(\sigma_{\rho},\sigma_n) > 0$.
\item
There exists a $\pi$ for which
$\liminf_{n \to \infty} \Delta_n(\pi) > \max_i \pi_i.$
\item
If $\pi$ is the uniform distribution on $\A$, then
$\liminf_{n \to \infty} \Delta_n(\pi) > 1/|\A|$.
\item
For all $\pi$ with $\min_i \pi_i > 0$, $(\sigma_n)_{n=1}^{\infty}$
has a non-trivial tail. 
\item
There exists a $\pi$ with $\min_i \pi_i > 0$ such that
$(\sigma_n)_{n=1}^{\infty}$ has a non-trivial tail.
\end{enumerate}
The analogous $8$ conditions are equivalent for $c_n$.
\end{Proposition}

\section{The Ising model} \label{sec:ising}
The only family of channels for which the problem is well understood
is the family of symmetric binary channels
\begin{equation} \label{eq:ising}
\M = \left( \begin{array}{ll} 1 - \eps & \eps \\
\eps & 1 - \eps \end{array} \right),
\end{equation}
where $\lam_2(M)$, the second largest (in absolute value) eigen value of $M$,
satisfies $\lam_2(M) = 1 - 2 \eps$.
Channel (\ref{eq:ising}) corresponds to the Ising model on the tree. 
The free measure for the Ising model on a finite tree is the probability
measure on configurations $\sigma$ of $\pm 1$, given by
\begin{equation} \label{eq:ising1}
\P[\sigma] = \frac{1}{Z} \exp(\sum_{v \sim w} \sigma_v \sigma_w),
\end{equation} 
where $Z$ is a normalizing constant.
The correspondence between (\ref{eq:ising}) and (\ref{eq:ising1}) is given
by
$\eps = \frac{\exp(-\beta)}{\exp(-\beta) + \exp(\beta)}$, or equivalently,
$\lam_2(M) = \tanh \beta$.

\begin{Theorem} \label{thm:ising}
The reconstruction problem is solvable for the binary
symmetric channel with error probability $\eps$ (\ref{eq:ising}), and the
$b$-ary tree $T_b$, if and only if $b \lam_2^2(M) > 1$.
If $b \lam_2^2(M) > 1$, then the reconstruction problem is also {\em census}
solvable.
\end{Theorem}

Reconstruction (and census) solvability when $b \lam_2^2(M) > 1$
was first proved in \cite{Hi} (\cite{KS} is earlier and does much more,
but is formulated in the language of multi-type branching processes).
\proofs
We think of the spin values as $\pm 1$.
Write $\lam = 1 - 2 \eps$, and let $S_n$ be the sum of the $\pm$ variables
at level $n$ of the tree.
Given that the spin at the root is $+$, the expected value of $S_n$ satisfies
\begin{equation} \label{eq:1st_moment}
\E^{+}[S_n] = \sum_{v \in L_n} \E^{+}[\sigma_v] = b^n \lam^n.
\end{equation}
Similarly, $\E^{-}[S_n] = -b^n \lam^n$.
The second moment of $S_n$ satisfies
\begin{equation} \label{eq:2nd_moment}
\E^{+}[S_n^2] = \E^{-}[S_n^2] = \sum_{v,w \in L_n} \E[\sigma_v \sigma_w]
              = b^n \left(1 + \sum_{j=1}^n (b^j - b^{j-1}) \lam^{2j} \right)
              = \Theta(b^{2n} \lam^{2n}),
\end{equation}
where the last equality follows from the fact that $b \lam^2 > 1$.
By Cauchy-Schwartz,
\begin{eqnarray*}
\E^{+}[S_n] - \E^{-}[S_n] 
 &=& \sum_{\sigma} (\P^{+}[\sigma] - \P^{-}[\sigma]) S_n(\sigma) \cr
 &\leq& 
  \sqrt{ \sum_{\sigma} \frac{(\P^{+}[\sigma] - \P^{-}[\sigma])^2}{\P^{+}[\sigma] + \P^{-}[\sigma]}}
  \sqrt{\sum_{\sigma} S_n^2(\sigma)(\P^{+}[\sigma] + \P^{-}[\sigma])}.
\end{eqnarray*}
It now follows by (\ref{eq:1st_moment}) and (\ref{eq:2nd_moment}) that when $b \lam^2 > 1$,
\[
\sum_{\sigma} \frac{(\P^{+}[\sigma] - \P^{-}[\sigma])^2}{\P^{+}[\sigma] + \P^{-}[\sigma]} = \Theta(1).
\]
which implies that $D_V(\P^{+},\P^{-}) = \Theta(1)$.
\QED

The reconstruction solvability result when $b |\lam_2(M)|^2 > 1$ is extended to general trees \cite{EKPS}
and general channels \cite{KS,MosselPeres:03},
where $b$ is replaced by the {\em branching number} of the tree and $\lam_2$ is the second largest eigenvalue
of the matrix $\M$ (in absolute value).

The proofs of the non-reconstruction result when $b \lam_2^2(M) \leq 1$ are harder, and do not generalize
to other channels. We know of 4 different proofs for this result
\begin{itemize}
\item
The first proof \cite{BRZ},
is based on recursive analysis of the Gibbs measure.
\item
A proof of non-reconstruction which is based on information inequalities
is given in \cite{EKPS} where it is shown that the mutual information
between the variable at the root of the tree and the level $n$ variables
satisfy $I(\sigma_{\rho},\sigma_n) \leq \sum_{v \in L_n} I(\sigma_{\rho},\sigma_v)$,
as in the case of conditionally independent variables.
This proof extends to general trees when $\br(T) \lam_2^2(M) < 1$,
where $\br(T)$ is the branching number of the tree.
\item
In \cite{I1} recursive analysis is used in order to show
$\E[\E^2[\sigma_{\rho} | \sigma_n]]$ tends to zero as $n \to \infty$ for the $n$-level
tree, when $b \lam_2^2(M) \leq 1$. The proof \cite{I2}
applies also to general trees when $\br(T) \lam_2^2(M) < 1$.
\item
Glauber dynamics is the following reversible Monte-Carlo
method for sampling configurations $\sigma$ according to the distribution 
(\ref{eq:ising}) or (\ref{eq:ising1}).
Given the current
configuration $\sigma$, a vertex $v$ is picked uniformly at random
at rate $1$,
in which case the variable $\sigma_v$ is replaced by a random variable $\sigma'_v$
chosen according to the conditional distribution on the
rest of the configuration, $(\sigma_w)_{w \neq v}$.
In \cite{BeKeMoPe:03} it is shown that Glauber dynamics have spectral
gap which is bounded away from zero when $b \lam_2^2(M) < 1$.
Then using a general principle (which is proven in a much more general context)
we obtain that the reconstruction problem is unsolvable when $b \lam_2^2(M) < 1$.
\end{itemize}
In \cite{PP} the critical case for general trees, $\br(T) \lam_2^2(M) = 1$,
is analyzed in detail.

\subsection{Reconstruction algorithms} \label{subsec:alg}
Theorem \ref{thm:ising} reveals a surprising phenomenon:
reconstruction by global majority vote has the same threshold for success
as maximum likelihood reconstruction (which is the optimal reconstruction strategy).

The parsimony method is popular in biology.
Given a bicoloring of the boundary of a tree $T$,
a {\em parsimonious} coloring of the internal nodes is any
assignment of the two colors to these nodes that minimizes the
total number of bicolored edges.
A way of finding a parsimonious coloring is the following:
Starting from the parents of the boundary
nodes, assign recursively
to each internal node the color of
the majority of its $\pm 1$-colored children.
In case of a tie, assign the non-color ``$?$''.
Then scan the tree from the root downwards and assign all vertices labeled by $?$ the same label as their
parent.

On a fixed finite tree, when $\eps > 0$ is small, the maximum likelihood
algorithm will reconstruct the same root value as one of the parsimonious
colorings given the boundary.

However, this is not the case when $\eps$ is larger.
For the binary tree, it is shown in \cite{St} that the parsimony
reconstruction algorithm has success probability bounded away from
$1/2$ as $n \to \infty$ if and only if $\eps \ge 1/8$.
Thus when $\lam_2(M) = 1 - 2 \eps\in (2^{-1/2},3/4]$ on the binary tree,
majority (and maximum likelihood) will have success probability bounded
away from $1/2$, while the parsimony success probability tends to $1/2$.

On a tree where each vertex has
$k$ children with $k$ odd, the above algorithm for finding a parsimonious coloring
reduces to recursive majority;
In~\cite{M:recur} it is shown that reconstruction via this method succeeds asymptotically if and only if
\[
\eps < \beta_k:= {1 \over 2} - {2^k \over 4k}
 {k-1 \choose {k-1 \over 2}}^{-1} \,.
\]
(it is interesting to note that the proof is based on exactly the same recursion which is analyzed in the context of noisy computation in \cite{HW,Ev}).

In~\cite{M:recur} more general reconstruction algorithms
on regular trees
(and more generally, on $\ell$-periodic trees) are analyzed.
Suppose that in order to determine
the color assigned to a node $v$, the algorithm is allowed
to examine the colors of its descendants $\ell$ generations down.
(However, only a single bit can be stored at each node). Then it is shown in \cite{M:recur} that
 recursively applying majority vote of the
descendants $\ell$ generations down is optimal,
yet it succeeds asymptotically only for flip probabilities $\eps$ below a threshold which is strictly lower
than the critical value for reconstruction.

Yuval Peres (private communication) conjectured that
\begin{Conjecture}
Consider the Ising model on the regular tree $T_b$ and reconstruction algorithms which are
{\em local:}
algorithms that for each vertex are allowed
to scan the information stored at its descendants $\ell$ generations down,
and at most $r$ bits of information
are allowed to be stored at each node.
Then for all $\ell$ and $r$
the threshold for reconstruction for such algorithms is strictly below the threshold
for reconstruction.
\end{Conjecture}
We remark that it is important to require that the algorithm examines only vertices below
the vertex which is being updated, as Glauber dynamics are local and are successful in
reconstruction of the root whenever $b \lam_2^2(M) > 1$.
A step of Glauber dynamics is performed as follows. 
Given the current
configuration $\sigma$, an internal vertex $v$ is picked uniformly at random
at rate $1$,
in which case the variable $\sigma_v$ is replaced by a random variable $\sigma'_v$
chosen according to the conditional distribution on the
rest of the configuration, $(\sigma_w)_{w \neq v}$.

We emphasize that the fact that recursive algorithms are asymptotically 
inferior to global majority does not hold for other Potts models or Ising models with
external fields(\cite{M:lam2}, see also Section \ref{sec:potts}).

\section{Census solvability} \label{sec:census}
The threshold $b \lam_2^2(M) = 1$ which appeared as the threshold both
for reconstruction solvability and for census solvability for the Ising model,
turns out to be in general the threshold for census solvability.
\begin{Theorem} \label{thm:blam2}
Let $M$ be a channel corresponding to an ergodic markov chain.
Let $T_b$ be the $b$-ary tree.
The reconstruction problem is census-solvable if $b |\lam_2(M)|^2 > 1$, and
is not census solvable if $b |\lam_2(M)|^2 < 1$.
For general trees, the reconstruction problem is solvable when $\br(T) |\lam_2(M)|^2 > 1$,
where $\br(T)$ is the branching number of the tree.
\end{Theorem}

\begin{Conjecture} \label{conj:census}
The reconstruction problem is not census solvable when 
$$b |\lam_2(M)|^2 = 1 \,.$$
\end{Conjecture}

{\bf $b |\lam_2(M)|^2 > 1$ implies census solvability.}
\cite{KS} proves a limit theorem for the variables $c_n$.
In particular it is shown, that if $b |\lam_2(M)|^2 > 1$ then the distribution of the
limiting variable depends on the initial variable at the root.
This implies that the problem is census solvable.

A more elementary proof is given in \cite{MosselPeres:03}.
The proof follows the lines of the proof for the
Ising model (Theorem \ref{thm:ising}), where $S_n$
is replaced by the scalar product of $c_n$ with any vector
nonzero $v$ satisfying $\M v = \lam_2(M) v$.
Note that this proof generalizes the proof for the Ising model,
since for the Ising model, $v = \left( \begin{array}{c} 1 \\ -1 \end{array} \right)$.
This proof also generalizes to general trees, proving that if $\br(T) |\lam_2(M)|^2 > 1$,
then it is possible to reconstruct the root using the scalar product of $v$ with a weighted
census $\tilde{c}_n(i) = \{\sum_x \omega(x) : x \in L_n, \sigma_x = i\}$ for some weights
$\{\omega(x)\}_{x \in T}$.

{\bf $b |\lam_2(M)|^2 < 1$ implies no census solvability.}
The CLT in \cite{KS} implies that if $b |\lam_2(M)|^2 \leq 1$ then the normalized
value of $c_n$ ($c_n/b^{n/2}$ if $b |\lam_2(M)|^2 < 1$)
converges to a nonzero random variable which is independent of
the variable of the root. However, this result on does not imply
that the reconstruction problem is not census solvable.
Presumably, it may the case that the first coordinate of $c_n$ is more
likely to be even for some value of the root variable than for others.
This dependency between the root variable and $c_n$ would not manifest itself in the
limiting normalized variables.

In \cite{MosselPeres:03} we combine the results of \cite{KS} with the local central limit
theorem to demonstrate that this could not happen.
The idea of the proof is to use \cite{KS} in order to couple
$c^i_n$, the value of $c_n$ given that the root variable is $i$, and
$c^j_n$, the value of $c_n$ give that the root variable is $j$
in such a way that the variables are close
(i.e., $|c^i_n - c^j_n|_{\infty} < \eps b^{n/2}$).
Then use the local central limit theorem in order to achieve a coupling
of $c^i_{n + \ell}$ and $c^j_{n + \ell}$ with high probability.
In \cite{MosselPeres:03} we also verify Conjecture \ref{conj:census} for Potts models
and asymmetric Ising models.

\section{Potts models} \label{sec:potts}
Two of the natural generalizations of binary symmetric channels are asymmetric
binary channels (which correspond
to Ising models with external field), and $q$-ary symmetric channels (which correspond to Potts models
with no external field):
\begin{itemize}
\item
Asymmetric binary channels have the state space
$\{0,1\}$ and the matrices:
\begin{equation} \label{eq:assymetric}
\M = \left( \begin{array}{ll} 1 - \delta_1 & \delta_1
                           \\ 1 - \delta_2 & \delta_2
           \end{array} \right),
\end{equation}
with $\lam_2(M) = \delta_2 - \delta_1$.
\item
Symmetric channels on $q$ symbols have the state space
$\{1,\ldots,q\}$ and the matrices:
\begin{equation} \label{eq:potts}
\M =
 \left(
\begin{array}{lllll} 1 - (q-1) \delta & \delta & \ldots     & \delta   \\
                    \delta & 1 - (q-1) \delta &  \delta       & \ldots \\
                    \vdots & \ldots & \ddots & \vdots            \\
                    \delta & \ldots & \delta     & 1 - (q-1) \delta \end{array}
        \right),
\end{equation}
with $\lam_2(M) = 1 - q \delta$.
\end{itemize}
Depending on the sign of $\lam_2(M)$ we distinguish
between {\em ferromagnetic}
Potts models where $\lam_2(M) > 0$,
and {\em anti-ferromagnetic} models where $\lam_2(M) < 0$.
When $1 - (q-1) \delta = 0$, we obtain the model
of proper colorings of the tree:
\begin{equation} \label{eq:color}
\M =
 \left(
\begin{array}{lllll} 0 & (q-1)^{-1} & (q-1)^{-1} & \ldots     & (q-1)^{-1}   \\
                    (q-1)^{-1} & 0 &  (q-1)^{-1}       & \ldots \\
                    \vdots & \ldots & \ddots & \vdots            \\
                    (q-1)^{-1} & \ldots & (q-1)^{-1}     & 0
\end{array}
        \right).
\end{equation}

\begin{Problem}
For the $3$ symbols Potts model (\ref{eq:potts}) find the
values for which the reconstruction problem is solvable
on the $b$-ary tree.
\end{Problem}

It may be easier to solve the analogous problem for
the Ising model with external field.
The analogous problem for colorings was stated in \cite{BW1}.
Applying standard coupon-collector
estimates recursively, it is
easy to see that if $b \geq (1+\delta) q \log q$ and $q$ is large, then the
reconstruction problem is solvable for the coloring model.

\begin{Problem}
For colorings, for which $b$ and $q$
is the reconstruction problem solvable on the $b$-ary tree?
\end{Problem}

Below we discuss several bounds for the reconstruction problem
for Potts models.
\begin{itemize}
\item
{\bf If $b \lam_2^2(M) > 1$ then the reconstruction problem is solvable.}
This follows from Theorem \ref{thm:blam2}, and from the fact that
census solvability implies solvability.
\item
{\bf If $b |\lam_2(M)| \leq 1$, then the reconstruction problem is
 unsolvable.}
\proofs
Assume first that $M$ is a ferromagnetic Potts model, i.e. $\lam_2(M) > 0$.
Consider two measures on the tree, one with $i$ as the root variable
and one with $j$ as the root variable. We {\em couple} these measures
in the following way: starting at the root if the two measures agree
on the variable at $v$, then we couple in such a way that the measures
also agree for all the children of $v$. If they do not agree at $v$, then for each of
the children of $v$, use the optimal coupling in order to couple the measures.
For each of the children, this has success probability $q \delta$. Thus the non-coupled vertices
are a branching process with parameter $1 - q \delta = \lam_2(M)$.
When $b \lam_2(M) \leq 1$ this process will eventually die; this means
that for large $n$ all the vertices at level $n$ will have the same variables
with probability going to $1$ as $n \to \infty$, as needed.
When $M$ is anti-ferromagnetic, the coupling probability is
$(q-2) \delta + 2(1 - (q-1) \delta) = 2 - q \delta$,
therefore the branching process parameter is
$1 - (2 - q \delta) = -\lam_2(M) = |\lam_2(M)|$.
Similar arguments apply for Ising models with external fields.
\QED
\item
{\bf If $b \lam_2(M) > 1$ and $q$ is sufficiently large, then the
reconstruction problem is solvable.}
This is the main result of \cite{M:lam2}. It implies in particular
that $b |\lam_2(M)|^2 = 1$ is not the threshold for the reconstruction
problem for Potts models, as it sometimes possible to reconstruct even when $b |\lam_2(M)|^2 < 1$.  
An analogous result is proven for the
asymmetric binary channel. The idea behind the proof is the following.
Channel (\ref{eq:potts}) may be thought of in the following way:
at each step the output is identical to the input with probability $\lam_2(M)$,
otherwise, the output is chosen uniformly among the $q$ symbols.
In particular if $\lam_2(M) > 0$ is fixed and $q$ is very large, then
if two of the children of a vertex in the $b$-ary tree $T_b$ have the same
label, then with overwhelming probability, this is also the label of their
parent. Now suppose that $q$ is large and there exists a copy of $T_2 \subset T_b$
such that all the vertices of $T_2$ are labeled by $i$.
Using a recursive argument we see that given this event, with large probability,
the variable at the root is $i$.
Moreover we show that when $b \lam_2(M) > C$ for some constant $C > 1$, such a tree exists with positive
probability.
Therefore, it is possible to reconstruct the root variable based on the existence of
such a unicolored tree.
In order to obtain the result for $C = 1$, we replace the unicolored $T_2$ by a diluted unicolored $T_2$.
\item
{\bf If $b \frac{(1 - q \delta)^2}{1 - (q-2) \delta} \leq 1$
 then the reconstruction problem is unsolvable.}
 This and the analogous result for asymmetric binary
 channels are proven in \cite{MosselPeres:03}, we sketch the main idea of the proof below.
\proofs
In order to show that the reconstruction problem is unsolvable, it suffices to show
that given that the root value is $0$ or $1$ with probability $1/2$ each,
it is asymptotically impossible to conclude from the variables at level $n$, if the root is $0$ or $1$.
Suppose that in addition to the variables at level $n$, we are also given all the
variables at all levels of the tree having variable $j$ with $j \neq \{0,1\}$.
Since we are given more information, it is easier to reconstruct.
The model where we are given this extra information is nothing but the symmetric binary
channel with matrix
\[
\M =
\left( \begin{array}{ll}
\frac{1 - (q-1) \delta}{1 - (q-2) \delta} & \frac{\delta}{1 - (q-2) \delta}\\
\frac{\delta}{1 - (q-2) \delta} & \frac{1 - (q-1) \delta}{1 - (q-2) \delta}
\end{array} \right),
\]
on the random tree which is obtained from the original tree by independently
deleting an edge with probability $(q - 2) \delta$ and retaining it with probability
$1 - (q - 2) \delta$.
The results of \cite{EKPS} imply that the binary symmetric channel on a general tree $T$ ,
the reconstruction problem is unsolvable if $\br(T) \lam_2^2(M) < 1$.
For the branching process on the regular tree $T_b$, the branching number is "typically"
$b \left(1 - (q - 2) \delta \right)$. The non-reconstruction criterion
$\br(T) \lam_2^2(M) < 1$, now translates to
$b \frac{(1 - q \delta)^2}{1 - (q-2) \delta} < 1$.
\QED
\end{itemize}
In Figure 1 we draw several of the bounds, where $b$ is a function of $\lam = \lam_2(M)$.
The area above the top curve, $b \lam^2 = 1$,
is the area where (census) reconstruction is successful for all channels.
The top curve is the critical curve for the symmetric binary channel:
above it reconstruction is successful and below it, it fails.
Below the second curve, $\lam = \frac{b (1 - 3 \delta)^2}{1 - \delta}$, reconstruction fails for
the $q=3$ Potts model.
Below the bottom curve $b \lam = 1$, reconstruction fails for all Potts models.
This curve is also the asymptotic critical curve as $q \to \infty$.
\begin{figure} [hptp] \label{fig:potts}
\epsfxsize=5cm
\hfill{\epsfbox{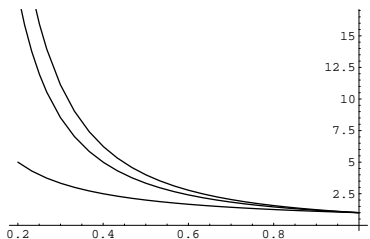}\hfill}
\caption{bounds for $b$ as a function of $\lam_2(M)$}
\end{figure}

\subsection{Algorithms}
We want to point out that the proofs in \cite{M:lam2} imply that for Potts models
when $q$ is large, reconstruction using recursive schemes has better threshold than
any algorithm which uses only the census.
We conjecture that the phenomenon that global majority achieves the same threshold for
reconstruction as maximum likelihood occurs only for symmetric binary channels.
\begin{Conjecture}
Consider Potts models (\ref{eq:potts}) for $q \geq 3$.
Then there exist $b$ and $\delta$ such that the reconstruction problem
is solvable for the $b$-ary tree, yet $b \lam_2^2(M) < 1$.
\end{Conjecture}
Note that this conjecture implies in particular, that any algorithm
which uses the census only is inferior to the optimal algorithm for Potts models.

\subsection{Monotonicity}
For Potts models (\ref{eq:potts}), it is easy to see that
if the reconstruction problem is solvable for $q$ and $\eps$
and $q' < q$, then the reconstruction problem is also solvable for $q'$
and $\eps$.
We expect that for fixed $\lam = \lam_2(M)$,
reconstruction is easier when $q$ is larger.
\begin{Conjecture}
Consider two symmetric channels $M_1$ and $M_2$
(as in  (\ref{eq:potts})) on $q_1$ and $q_2$
symbols respectively, where $q_1 < q_2$.
If $0 < \lam_2(M_1) = \lam_2(M_2)$ and the reconstruction problem
 is solvable for $M_1$, then it is also solvable for $M_2$.
\end{Conjecture}
This is obvious when $q_2$ is a multiple of
$q_1$. Using the reconstruction criterion for the binary symmetric channel on
$2$ symbols, it is easy to prove the conjecture when $q_1 = 2$.

\section{General channels} \label{sec:strange}
In this section we discuss techniques which apply to the reconstruction
problem for general channels.
\subsection{Proving solvability}
{\bf Spectral methods.} Theorem \ref{thm:blam2} implies that
when $b \lam_2^2(M) > 1$ the reconstruction problem is census-solvable
(and therefore solvable) for the $b$-ary tree and the channel $M$.

\medskip\noindent
{\bf Recursive methods.} Starting at the boundary of the tree,
we may try to evaluate recursively the variable at each vertex
of the tree. Assuming that the probability of reconstructing the
value of a variable at distance $n$ from the boundary is $p_n$, we
obtain recursive bounds on $p_{n+1}$. Bounding these recursions,
we may prove that the reconstruction problem is solvable.
Although for symmetric binary channels, these methods
always achieve worse thresholds than the spectral method \cite{M:recur},
in many other cases the recursive methods are superior.
The first example which was already discussed in Section
\ref{sec:potts} is that of the Potts model on $q$ symbols where $q$ is
large. It is proven in \cite{M:lam2} that if $b \lam_2(M) > 1$ and
$q$ is sufficiently large then the reconstruction problem is solvable
(the spectral method only applies when $b \lam_2^2(M) > 1$).
In \cite{MosselPeres:03} we note that a similar argument proves that for the coloring
problem for large $q$, if $b \geq (1 + \eps) q \log q$, then the reconstruction
problem is solvable.
In general this method combined with standard large deviation estimates (see e.g. \cite{DZ})
implies:
\begin{Theorem} \label{thm:gen}
Let $M$ be a channel, such that for all $i,j \in \A$, there exists
an $\ell$ such that $\M_{i,\ell} \neq \M_{j,\ell}$.
Then when $b$ is sufficiently large, the reconstruction problem is
solvable for the tree $T_b$ and $M$.
\end{Theorem}
This is proven in \cite{MosselPeres:03} where a general criterion is given to decide:
does there exist for the channel $M$ a number $b$ such that the reconstruction
problem is solvable for $T_b$ and $M$.

\begin{Example} \label{ex:shift}
Let $\{Z_i\}_{i=1}^{\infty}$ be an i.i.d. sequence of variables such that
$\P[Z_i = 0] = 1 - P[Z_i = 1] = p$ where $0 < p < 1$.
Let $h \geq 1$ and consider the channel $M$ defined by the markov chain
$Y_i = (Z_i,\ldots,Z_{i+h})$.
Thus $M$ has state space $\{0,1\}^{h+1}$ with the product $(p,1-p)$
probability measure.
It is easy  to see that for the tree process, the variable $(\sigma_v)_{|v| \leq n}$ are independent
of $(\sigma_v)_{|v| \geq n+h}$, and therefore the reconstruction problem is not
solvable for $M$.
On the other hand, letting
\[
Y_i = \max_{0 \leq j \leq h} \{Z_i = \cdots = Z_{i+j} = 1\}.
\]
It is easily seen that $Y_i$ defines a channel $M$ on the space
$\{0,\ldots,h\}$.
Moreover, it is clear that for all $\ell$ the variables $\{Z_i\}_{i \geq \ell+h+1}$ and $\{Y_i\}_{i \leq \ell}$
are independent.
Therefore $\{Y_i\}_{i \geq \ell+h+1}$ and $\{Y_i\}_{i \leq \ell}$ are independent.
It follows that the variables $M^{h+1}(j)$ have the same distribution for all $j$.
Thus $\rank(\M^{h+1}) = 1$, and $\lam_2(M) = 0$.
Writing $\M$:
\[
\M = \left( \begin{array}{lllll}    p &   p(1-p) &    p (1-p)^2 &     \ldots &    (1-p)^h \\
                                    1   &   0 &         \ldots &             \ldots &    0   \\
                                    0   &   1   &       0 &             \ldots &         \vdots \\
                                    \ldots & \ddots & \ddots &         \ddots &    \vdots \\
                                    0 &     \ldots  &  0     &          p &       1-p \\
            \end{array} \right),
\]
if follows from Theorem \ref{thm:gen} that the reconstruction problem is solvable
for $M$ and $T_b$ provided that $b$ is sufficiently large.
This is a generalization of a channel appearing in \cite{M:recur};
see also \cite{LP2}.
\end{Example}
The example above demonstrates that it may be the case that for
the markov chain corresponding to the channel $M$, the states at times $t$ and $t+h$
are independent, yet, for the tree process, the reconstruction problem is solvable.
In fact, a much stronger phenomenon occurs
\begin{Theorem} \label{thm:crypt}
Let $b > 1$ be an integer and $T$ be the $2$-level $b$-ary tree.
There exists a channel $M$ such that for any initial distribution,
$\sigma_\rho$ and $\sigma_{\partial}$ are independent
(where $\sigma_{\rho}$ is the root label, and $\sigma_{\partial}$ is the
configuration at the leaves of the $2$-level $b$-ary tree), yet
when $B$ is sufficiently large, the reconstruction problem for the channel $M$
and the infinite $B$-ary tree $T_B$ is solvable.
\end{Theorem}
The construction in \cite{MosselPeres:03} is motivated by work on
secret-sharing protocols \cite{Sh} and applies Theorem \ref{thm:gen}.
We define the channel below. For the proof we refer the reader to
\cite{MosselPeres:03}.

\medskip\noindent
{\bf construction.}
Let $\F$ be a finite field with $q > b+2$ elements.
Let $x_1,\ldots,x_{b+1}$ be a fixed set of non-zero elements of $\F$.
We define a channel on the state space
\[
\F^b[x] = \{ f(x) : f(x) \in \F[x], \deg f \leq b \}.
\]
Given $f$, take $I$ to be a uniform variable in the set $\{1,\ldots,b+1\}$,
then take $M(f)$ to be $g \in \F^b[x]$ chosen 
uniformly among the $g$'s satisfying $g(0) = f(x_I)$.
Given the value of $f$ at $b$ of the points $x_1,\ldots,x_{b+1}$, and for all 
$a \in \F$, there exists a unique polynomial satisfying $f(0) = a$.
This implies that $b$ independent copies of the chain at $f$ give no 
information on $f$. In \cite{MosselPeres:03} it is shown that when the number of copies
is sufficiently large, information is retained so that the reconstruction 
problem is solvable.  

\subsection{Proving non-solvability}
We have a few techniques for proving non-solvability.

\medskip\noindent
{\bf Spectral gap of Glauber dynamics.}
Let $\Lambda(n)$ be the spectral gap of Glauber dynamics for the $n$-level $b$-ary tree. 
In \cite{BeKeMoPe:03} we prove a result for general graphs which implies for $T_b$ the following:
\begin{Theorem}
Suppose that $M$ is a channel such that Glauber dynamics satisfy
$\inf_n\Lambda(n) > 0$, then the mutual information between
$\sigma_{\rho}$ the root variable of $T_b$, and $\sigma_n$,
the variables at level $n$ decays exponentially fast: $I(\sigma_{\rho},\sigma_n) = O(\exp(-\Omega(n)))$.
In particular, the reconstruction problem for the tree $T_b$ is unsolvable.
\end{Theorem}
\begin{Problem}
Does it hold for reversible $M$ and the $b$-ary tree $T_b$
that $I(\sigma_{\rho},\sigma_n) = O(\exp(-\Omega(n)))$ if and only if
$\inf_n \Lambda(n) > 0$?
\end{Problem}
For the Ising model on trees, the answer to the problem is positive,
see \cite{BeKeMoPe:03}.

\medskip\noindent
{\bf Recursive analysis of maximum likelihood.}
A direct approach to the reconstruction problem is
to analyze the distribution of the ($\log$) likelihood
of the root variable given the boundary variable.
This leads to an iteration of random variables.
The only case in which this iteration was analyzed
is the symmetric binary channel (\cite{PP}) where this approach
yields an exact criterion for reconstruction for general trees.
It is an interesting challenge to extended
this technique to other channels.

It may be easier to analyze these recursions for "robust" phase transitions which
first appeared in \cite{PS}.
Consider the usual reconstruction problem, but suppose that the data at the
boundary is given with some additional noise.
The proofs that if $b \lam_2^2(M) > 1$ the reconstruction problem is
(census) solvable are immune to this noise.
However, this may not be the case for the reconstruction problem.
Indeed, we suspect that adding this additional noise
(assuming it is fixed but sufficiently strong) will shift the
phase transition to the point $b \lam_2^2(M) = 1$.
A similar phenomenon was proven in \cite{PS} for the phase transition
of uniqueness.
For $n$ and $m$, we denote by $\sigma_{n,m}$ the configuration which
is obtained from $\sigma_n$ by applying the random function
$M^m$ independently
on each of the symbols in $\sigma_n$. We denote by $\P^{\ell}_{n,m}$
the conditional distribution of $\sigma_{n,m}$ given that
$\sigma_{\rho} = \ell$. We then
\begin{Conjecture}
For all $M$ and $b$, such that $b \lam_2^2(M) < 1$, there exists $m$
such that for the $b$-ary tree
\[
\sup_{i,j} \lim_{n \to \infty} D_V(\P^i_{n,m},\P^j_{n,m}) = 0.
\]
\end{Conjecture}

\section{Terminology and related problems} \label{sec:terminology}
\subsection{Related problems}
In this subsection several variants of the reconstruction problem
are discussed.
Throughout the section we will
assume that the variable at the root is chosen uniformly.
By Proposition \ref{prop:equiv},
reconstruction solvability is equivalent
to the fact that there exists $\delta > 0$ such that for all $n$ with {\em probability
at least $\delta$}, the conditional distribution of $\sigma_{\rho}$ given $\sigma_n$
has {\em total variation distance at least $\delta$} from the uniform distribution.

We may consider the following variants of the problem:
\begin{itemize}
\item
{\bf non-uniqueness of the Gibbs measure.}
The condition here is that for all $n$ there {\em exists} $\sigma_n$ such that the
distribution of $\sigma_{\rho}$ given $\sigma_n$ has {\em total variation distance at least
$\delta > 0$} from uniform. This is a weaker condition than reconstruction solvability and it
was studied in statistical physics for Ising and Potts models. In particular,
the phase transition for these models is known, see \cite{Ha}.

\item
{\bf dismantlable graphs.} Suppose that we require that for all $i \in A$ and for all $n$
there {\em exists} $\sigma_n$ such that {\em $\P[\sigma_{\rho} = i | \sigma_n] = 1$}.
This requirement clearly
fails for all trees if the matrix $\M$ satisfies $\M_{i,j} > 0$ for all $i$ and $j$.
Moreover, this property depends only on which of
the entries of $\M$ are non-zero. Define a directed graph $G$ on $\A$ such that
$(i,j)$ is an edge of $G$ {\bf iff} $\M_{i,j} > 0$. We claim that there exists for sufficiently large $b$ a $b$-ary tree $T$
for which the requirement holds {\bf iff} for all $i \neq j$ the sets $N(i) = \{\ell : (i,\ell) \in G\}$
and $N(j) = \{\ell : (j, \ell) \in G\}$ satisfy
\begin{equation} \label{eq:Nrelation}
N(i) \not\subset N(j).
\end{equation}

{\bf proof.}
If $N(i) \subset N(j)$ then for all $n \geq 1$ there exists no $\sigma_n$ such that
given $\sigma_n$ the value of $\sigma_{\rho}$ is $i$ with probability $1$.
On the other hand, if for all $i \neq j$, it holds that $N(i) \not \subset N(j)$,
then given $i$, consider the labeling $\sigma_n$ of $T_b$ for $b \geq |\A|$,
which is obtained in the following way:
The root satisfies $\sigma_{\rho} = i$. Given $\sigma_v$, label the children of
$v$, denoted $w_1,\ldots,w_b$, in such a way that
$\{ \sigma_{w_i} : 1 \leq i \leq b \} = \{\ell , (\sigma_v,\ell) \in G\}$.
For all $n$ it now holds that given $\sigma_n$ the value of the root $\sigma_{\rho}$ must be $i$.
\QED

If (\ref{eq:Nrelation}) holds for all $i \neq j$,
then for large $b$ for $\delta > 0$
fraction of the $\sigma_n$,
it holds that $\P[\sigma_{\rho} = i | \sigma_n] = 1$
for some $i$ which depends on $\sigma_n$.
This may be proved using a recursive argument similar
to \cite[Theorem 2.1]{MosselPeres:03}.

We may replace the above requirement by the requirement
that there exists $i \in A$ such that for all $n$
there exists $\sigma_n$ such that given $\sigma_n$
the variable $\sigma_{\rho}$ satisfies
$\sigma_{\rho} \neq i$ (with probability $1$).
Defining $G$ as in the previous case, the property holds for $T_b$ for
sufficiently large $b$ {\bf iff} $G$ is {\em dismantlable}, see \cite{BW2}.

\item
{\bf census.} One may ask similar questions about the census.
\begin{itemize}
\item
Does there exists a $b$, such that for the tree $T_b$,
for all $n$ there {\em exists a census $c_n$} such that the distribution
of $\sigma_{\rho}$ given $c_n$ has {\em total variation distance at least $\delta$
from uniform}? For some models (like ferromagnetic Ising and Potts models) this
condition is equivalent to uniqueness of the Gibbs measure. For others (like colorings)
it seems that these conditions are not equivalent (we do not know how to demonstrate it
for colorings; Example \ref{ex:shift} is an example of such model).
\item
Does there exists a $b$, such that for the tree $T_b$
and all $n$ there exists a census $c_n$ such that
$\P[\sigma_{\rho} \neq i | c_n] = 1$ for some $i$?
It follows from \cite[Lemma 6.2]{MosselPeres:03} that such a $b$ does not exist
when $M$ is ergodic.
\end{itemize}
\item
{\bf phylogeny.} Most of the biological research which is related to 
reconstruction is devoted to problems in which the underlying tree is unknown and the algorithm is supposed to
find both the tree and the variables at the nodes of the tree, see e.g. \cite{Ca,Fi,Ha}.
These problems seems to be quite hard to analyze; in particular, in some cases there is no
well-defined probability space of trees.
In a recent work \cite{Mossel:03a,Mossel:03b} 
we show that the reconstruction of
phylogenetic trees is closely related to the reconstruction problem. 
 
\item
{\bf Noisy computation.}
Von Neumann~\cite{VN} proposed a model of computation in noisy
circuits where each gate computes correctly with probability $1-\eps$,
The analysis of this model in \cite{VN,Pi,Ev,ES1} has many similarities
to the analysis of the reconstruction problem for the symmetric binary channel,
in \cite{BRZ,I1,I2,EKPS,M:recur}.
However, we do not know of any formal relationship between the two models.
\end{itemize}

\subsection{Dictionary}
As the reconstruction problem has been studied from different perspectives,
different terminology is often used for the same entities.
We list below some terms and their translations.

\medskip\noindent
{\bf variable.} also symbol, state (finite markov chains), label, letter,
message (information theory), spin (statistical physics), color (combinatorics),
genotype (biology), phenotype (biology).

\medskip\noindent
{\bf $b$-ary trees.} also $b+1$ regular trees (corresponding to the degree as a graph), and Bethe lattice (statistical physics).

\medskip\noindent
{\bf channel.} the channel (random function) $M$ corresponds to a stochastic
matrix $\M$ such that $\M_{i,j} = \P[M(i) = j]$. In the statistical physics literature,
when working with Ising and Potts models,
it is common to work with the {\em Hamiltonian}:
\begin{equation} \label{eq:hamiltonian}
H((\sigma_v)_{v \in T}) = \sum_{v} h_{\sigma_v} +
\beta \sum_{(v,w) \mbox{ edge of } T} \delta_{\{\sigma_v = \sigma_w\}}.
\end{equation}
The probability of a configuration $\{\sigma_v\}_{v \in T}$ is then
\begin{equation} \label{eq:partition}
\frac{1}{Z} \exp \left ( H((\sigma_v)_{v \in T}) \right),
\end{equation}
where $Z$ is a normalizing constant, known as the {\em partition function} (it is a function
of $H$ and of the tree $T$).
The parameter $1/\beta$ is often referred to as the temperature.

From (\ref{eq:hamiltonian}) the matrix $\M$ is given by:
\[
\M_{i,j} = \frac{\exp \left( h_j + \beta \delta_{\{i = j\}} \right)}
{\sum_{\ell} \exp \left( h_{\ell} + \beta \delta_{\{i = \ell\}} \right)}
\]
Some families of interest are:
The Ising model with no external field, when
$|\A| = 2$ and $h_1 = h_2 = 0$;
The Ising model with external field, where
$|\A| = 2$; and Potts models without external field where $|\A| = q$
and all the $h$ values are $0$.

The process on the infinite tree then corresponds to the {\em Gibbs} measure
on that tree with the specification that the root distribution is uniform.

In biology, the matrix $\M$, is related to the {\em mutation rate} of the process, and it
is usually assumed that $\M$ is a perturbation of the identity matrix.

In combinatorics the popular model is proper colorings of the tree, where
$\M$ is a $q \times q$ matrix satisfying $\M_{i,j} = \delta_{\{i \neq j\}} (q-1)^{-1}$.
This model is also referred to as the zero temperature anti-ferromagnetic Potts model,
since $\M$ is obtained as the limit of the corresponding $\M$ for Potts models when
$\beta \to -\infty$.

\medskip\noindent
{\bf reconstruction solvability.} Corresponds in statistical physics 
to extremality of the above measure. 
In information theory, it is natural to express reconstruction solvability
in term of decay of mutual information (see Proposition \ref{prop:equiv}).

\section{Very recent results}
Since the submitting this survey, a number of new results on 
reconstruction appeared. In \cite{MaSiWe:03b} better bounds for the 
reconstruction problem for Potts and asymmetric binary channels were obtained.
Essentially the same bounds for asymmetric binary channels were obtained 
independently in \cite{Martin:03}.  
\cite{MaSiWe:03b} provides a comprehensive 
analysis of the mixing rates of Glauber dynamics for Ising and Potts models 
under various boundary conditions and gives a positive answer to Problem 3.  

In \cite{JansonMossel:04} we analyze robust solvability. 
This may be thought of as reconstruction where the labels at the bottom level
are further perturbed. We show that the threshold for robust reconstruction 
is given by $b |\lam_2(M)|^2 = 1$ as conjectured in \cite{MosselPeres:03}. 

Finally, we note that the crucial role of the reconstruction 
problem in Phylogeny was recently demonstrated in \cite{Mossel:03b} 
and \cite{MosselSteel:04} extending the results of \cite{Mossel:03a}.

\medskip\noindent
{\bf Acknowledgments:}
I learned about the reconstruction problem from Yuval Peres.
I want to thank him for many fruitful conversations about the problem.
I thank Olle H\"{a}ggstr\"{o}m, Claire Kenyon,
L\`aszl\`o Lov\`asz, Jeff Steif and Peter Winkler for helpful discussions, and
the referee for many helpful comments.



\begin{thebibliography}{7}
\bibitem{I1} Ioffe, D. (1996a).  A note on the extremality of the disordered
state for the Ising model on the Bethe lattice.
{\em Lett.\  Math.\  Phys.} {\bf 37}, 137--143.
\end{thebibliography}
\end{document}